\documentclass[12pt]{article}  

\renewcommand{\baselinestretch}{1.3}
\normalsize

\newcommand{\Z}{{\bf Z}}

\newcommand{\C}{{\bf C}}
\renewcommand{\S}{{\bf S}}
\renewcommand{\P}{{\bf P}}

\usepackage{amsmath}
\usepackage{epsfig}
\usepackage[psamsfonts]{amsfonts}

\newcommand{\g}{\mathfrak g}
\newcommand{\h}{\mathfrak h}

\newcommand{\be}{\begin{equation}}
\newcommand{\ee}{\end{equation}}
\newcommand{\bea}{\begin{eqnarray}}
\newcommand{\eea}{\end{eqnarray}}
\newcommand{\no}{\noindent}

\newcommand{\G}{\Gamma }
\renewcommand{\L}{\mathfrak{L}}
\newcommand{\bp}{\begin{para}}
\newcommand{\ep}{\end{para}}
\newcommand{\su}{\mathfrak{su}(2)}
\newcommand{\bd}{\begin{defin}}
\newcommand{\ed}{\end{defin}}
\newcommand{\cz}{\C^3/\Z_3}
\newcommand{\CZ}{\C^n/\Z_n}
\newcommand{\non}{\nonumber}
\newcommand{\ggg}{\g_\L}
\newcommand{\sm}{SU(3)\times SU(2)\times U(1)}

\numberwithin{equation}{section}

\begin{document} 
\newtheorem{lfact}{Little Fact}[section]

\newtheorem{theorem}{Theorem}[section]
\newtheorem{proposition}{Proposition}
\newtheorem{conjecture}{Conjecture}

\newcounter{def}[section]

\newcounter{defin}[section]
\newenvironment{defin}
{\vskip .3cm \noindent \stepcounter{defin}{\bf Definition} 
{\bf {\arabic{section}.{\arabic{defin}} }}
\sl 
}

\newenvironment{fact}
{\vskip .3cm \noindent \stepcounter{def}{\bf Fact} 
{\bf {\arabic{section}.{\arabic{def}} }}
\sl 
}

\newcounter{eg}[section]
\newenvironment{eg}
{\vskip .3cm \noindent \stepcounter{def}{\bf Example }
{\bf {\arabic{section}.{\arabic{def} }}}}

\newcounter{lemma}[section]

\newenvironment{lemma}
{\vskip .3cm \noindent \stepcounter{def}{\bf Lemma} 
{\bf {\arabic{section}.{\arabic{def}} }}
\sl 
}

\newenvironment{para}[1]
{\vskip .3cm {\bf #1.}\hskip .4cm}

\newcommand{\qed}{\nobreak \ifvmode \relax \else
      \ifdim\lastskip<1.5em \hskip-\lastskip
      \hskip1.5em plus0em minus0.5em \fi \nobreak
      \vrule height0.75em width0.5em depth0.25em\fi}

\renewcommand{\baselinestretch}{1.5}
\normalsize

\thispagestyle{empty}

\Large

\begin{center} 

{\bf {Orbifold Singularities, 

Lie Algebras of the Third Kind (LATKes),

and Pure Yang-Mills with Matter}}

\end{center}
\vskip 1cm
\large
\begin{center} {\tt Tamar Friedmann}\footnote{E-mail: tamarf@mit.edu}$^,$\footnote{Current address: tamarf@pas.rochester.edu}\end{center}
\normalsize 
\centerline{\it Massachusetts Institute of Technology}
\centerline{\it Cambridge, MA 02139, USA}

\vskip 1.5cm

\abstract{
We discover the unique, simple Lie Algebra of the Third Kind, or LATKe, that stems from codimension 6 orbifold singularities and gives rise to a  new kind of Yang-Mills theory which simultaneously is pure and contains matter. 
The root space of the LATKe is 1-dimensional and its Dynkin diagram consists of one point.
The uniqueness of the LATKe is a vacuum selection mechanism. 
}

\vskip 1cm

\thispagestyle{empty}
\newpage

\begin{center}
\Large

The World in a Point ?
\end{center}
\begin{quote}\end{quote}

\vskip 4.5cm

\renewcommand{\baselinestretch}{1.2}
\normalsize


\begin{center}
\hskip 2cm Blow-up  of $\cz$ $|$  Dynkin diagram  of the LATKe \hskip 1cm
\end{center}

\vskip -0.2cm
\begin{center}
\Huge
\[ \bullet \]
\end{center}

\normalsize

\vskip 2.1cm

\begin{center}
Pure Yang-Mills 
with matter

\end{center}
\normalsize


\renewcommand{\baselinestretch}{1.4}
\normalsize

\thispagestyle{empty}

\newpage

\setcounter{page}{1}

\section{Introduction}

The interpretation of codimension 4  orbifold singularities  as ADE gauge theories, which arose within string theory in the mid 90's \cite{Hull:1995mz,Witten:1995ex}, has been extended to the case of M--theory compactifications, where codimension 4 orbifold singularities in $G_2$ spaces were also interpreted as ADE gauge theories  \cite{Acharya:1998pm,Acharya:2000gb,Atiyah:2001qf}. Further orbifolding the already-singular $G_2$ spaces led us to the first manifestation via M--theory of Georgi-Glashow grand unification:  from an $A_4$ singularity of the $G_2$ space, an $SU(5)$ gauge theory broken by Wilson lines  precisely to the gauge group of the standard model $\sm$ arose naturally, with no extraneous gauge fields \cite{Georgi:1974sy,Friedmann:2002ct,Friedmann:2003cd}. A precise relation between the energy scale of grand unification ($M_{GUT}$) and certain volumes inside the $G_2$ space was also obtained \cite{Friedmann:2003cd,Friedmann:2002ty}.

In the process of constructing the $G_2$ spaces, orbifold singularities of codimension 6 arose as well \cite{Friedmann:2002ct}. However, there was no analog of the interpretation of a codimension 4 singularity as an ADE gauge theory for the case of codimension 6. We set out to find such an analog.

To our delight and surprise, we discover far more than we  expected, both mathematically speaking and physically speaking. 

On the math side, we introduce a new set of relations, which we call the Commutator-Intersection Relations, that illuminate 
the connection between codimension 4 singularities and Lie algebras. These relations pave the way to construct Lie Algebras of the Third Kind, or LATKes, a 
kind of algebras that arise from codimension 6 orbifold singularities. We also learn and prove the existence and uniqueness of  a simple LATKe.

On the physics side, we discover a new
kind of Yang-Mills theory, called "LATKe Yang-Mills," which arises from the LATKe.  Unlike any known Yang-Mills theory, the LATKe Yang-Mills theory in its purest form automatically contains matter. We also propose that the uniqueness of the simple LATKe is a vacuum selection mechanism. The selected vacuum theory is  an $SU(2)\times SU(2)$  gauge theory with matter in the $(2,2)$ representation, and the  corresponding singularity is $\cz$. The algebra $\su \times \su$ is protected by the LATKe from being broken. The selected singularity  $\cz$ is one of those which arose in the $G_2$ spaces of \cite{Friedmann:2002ct}, and which at the time we  put on hold in anticipation of  the outcome of this investigation.

\section{The Codimension 4 Case}\label{codim4}

In this section, we review the correspondence between codimension 4 orbifold singularities and ADE Lie algebras, introduce the Commutator-Intersection Relations, and review the relation between physical interactions on the one hand and commutators and intersections on the other hand.

\subsection{Du Val-Artin correspondence}\label{DuVA}
The interpretation of codimension 4 orbifold singularities as ADE gauge theories is mathematically rooted in the work of Du Val and of Artin \cite{duval,artin,steinberg}, who  pointed out a correspondence between certain singularities and their blow-ups on the one hand and certain Lie algebras on the other hand. 

Before we state the correspondence, we provide below the necessary ingredients. 

The singularities in question are those that appear at the origin of $\C ^2$ under the orbifold action of finite discrete subgroups of $ SU(2)$. 
These subgroups, denoted $\G$,  had been classified as early as 1884 by F. Klein \cite{klein}. They consist of the cyclic groups $\Z _n$, also denoted $A_{n-1}$; the binary dihedral groups ${\bf D}_n$; and three "exceptional" groups: the binary tetrahedral ${\bf T}$, binary octahedral ${\bf O}$, and binary icosahedral ${\bf I}$, also denoted $E_6, E_7$, and $E_8$, respectively. Such a classification is known as an ADE classification. 

Each of these subgroups of $SU(2)$ has a natural action on $\C ^2$. For example, $\Z_n$ is generated by the $SU(2)$ matrix 
\be \label{Zn} \left ( \begin{array} [h] {cc} e^{2\pi i/n}&0 \\ 0&e^{-2\pi i/n} \end{array}\right )
\ee
and acts on $(x, y)\in \C ^2$ via the two--dimensional representation
\be \label{2Drep} (x, y)\longmapsto (e^{2\pi i/n}x, e^{-2\pi i/n}y)~.
\ee

The singularity at the origin of $\C ^2$ is analyzed by blowing up:  the singular space is replaced by a smooth manifold that looks just like $\C^2/\G$ everywhere except at the origin, and the origin itself is replaced by a smooth space of real dimension 2. This 2-dimensional space, known as the exceptional divisor, turns out to be a union of intersecting $\P ^1$'s, or 2--spheres $\S ^2$.

In the ambient four-dimensional space, the intersection of any two $\P ^1$'s is zero dimensional, i.e. it is a set of points. Counting those points gives an intersection number. The set of all intersection numbers forms the intersection matrix of the exceptional divisor, which we denote $\{ I_{ij}\}$. The indices $ij$ run from 1 to $b_2$, where $b_2$ is the number of independent 2--cycles in the exceptional divisor.

As it happens, the intersection matrix of the exceptional divisor of $\C^2/\G$ is equal to the negative of the Cartan matrix of the A, D, or E Lie algebra corresponding to $\G$:
\be \label{cequalsi} C_{ij}=-I_{ij} ~.\ee
In addition, the exceptional divisor itself becomes the Dynkin diagram of the corresponding Lie algebra when we replace each component $\P ^1$ of the exceptional divisor by a node and connect a pair of nodes by an edge when their corresponding $\P ^1$'s intersect; the components of the exceptional divisor therefore correspond to the simple positive roots of the Lie algebra. For example, when $\G = \Z _3$, the exceptional divisor is two intersecting $\P ^1$'s and the Dynkin diagram consists of two connected nodes:
\be \bullet \frac{ \mbox{  \hskip 1cm }  }{\mbox{  \hskip 1cm            }   } \bullet ~.
\ee

Now we are ready to state the Du Val-Artin correspondence: 

\bea \label{reptop}
\left \{
\begin{array}{rr} \label{divDynk}
\text {Exceptional divisor}\\
\text{of blow-up of $\C ^2/\G$} 
\end{array}
\right \}
&=&
\left \{
\begin{array}{ll}
\text {Dynkin diagram of}\\
\text {ADE Lie algebra}
\end{array}
\right \}
\\ 
\nonumber &&\\ 
 \left \{
\begin{array}{rrc} \label{intCart}
\text {Intersection matrix of }\\
\text{exceptional divisor} \\
I_{ij} \hskip 1cm 
\end{array}
\right \}
&=&
\left \{
\begin{array}{llc}
\text {Negative Cartan matrix}\\
\text {of ADE Lie algebra}\\
\hskip 1cm -C_{ij} 
\end{array}
\right \}
\\ \nonumber &&
\eea

\subsection{Commutator-Intersection Relations}\label{CIsection}
Here we rephrase the Du Val-Artin correspondence  as a set of relations between commutators of the Lie algebra and intersection numbers of the exceptional divisor, as follows. 

A complex simple Lie algebra is generated by $k$ triples $\{X_i, Y_i, H_i \}_{i=1}^k$ with their commutators determined by the following relations:

\bea \label{Serre}&&[H_i, H_j]=0\; ; \non 
\\&& [X_i, Y_j] = \delta _{ij} H_j \; ;\non
\\
&&[H_i, X_j]=C_{ij}X_j \; ; 
\\&&  [H_i, Y_j]=-C_{ij}Y_j \; ;\non
\\
&&\mbox{ad}(X_i)^{1-C_{ij}}(X_j)=0\; ; \non
\\ 
&&\mbox{ad}(Y_i)^{1-C_{ij}}(Y_j)=0\; .\non
\eea
Here, the $H_i$ form the Cartan subalgebra, the $X_i$ are simple positive roots, the $Y_i$ are simple negative roots, $k$ is the rank of the Lie algebra, $C_{ij}$ is the Cartan matrix, and  ad$(X_i)(A)=[X_i, A]$. 
These equations are the familiar Chevalley-Serre relations.

Now recall from equations (\ref{cequalsi}) and (\ref{intCart}) that $C_{ij}=-I_{ij}$. If we replace  $C_{ij}$ in equations (\ref{Serre}) by $-I_{ij}$, we get a new set of relations:
\bea \label{CI}&&[H_i, H_j]=0\; ; \non 
\\&& [X_i, Y_j] = \delta _{ij} H_j \; ;\non
\\
&&[H_i, X_j]=-I_{ij}X_j \; ; 
\\&&  [H_i, Y_j]=I_{ij}Y_j \; ;\non
\\
&&\mbox{ad}(X_i)^{1+I_{ij}}(X_j)=0\; ; \non
\\ 
&&\mbox{ad}(Y_i)^{1+I_{ij}}(Y_j)=0\; .\non
\eea

These relations demonstrate that the intersection numbers of the exceptional divisor  
completely determine the commutators of the corresponding Lie algebra.

This role of the intersection numbers in the structure of the Lie algebra will be central for us in what follows, and we will refer to  the relations (\ref{CI}) as the {\sl Commutator-Intersection Relations}, or the {\sl CI Relations}.

\subsection{Interactions, commutators, and intersections}\label{GE}

Physically speaking, there are relations between interactions and commutators, and interactions and intersections.

Here, we explain  roughly  how ADE gauge fields  arise from the codimension 4 singularities \cite{Mayr:1998tx,Katz:1996fh,vafa}.
First,  Kaluza Klein reduction of 3-form C--fields on 2--cycles of the exceptional divisor provides the gauge fields corresponding to the Cartan subalgebra.
Second, D2--branes wrapped on 2--cycles provide the "charged" gauge fields, forming the rest of the Lie algebra. 

For example, in the $\C^2/\Z_2$ case, where the exceptional divisor is a single $\P ^1$, three fields arise: a 3--form field reduced on the $\P^1$, which gives rise to the Cartan element denoted $A_\mu$; and D2--branes wrapped on the $\P^1$, which give rise to two oppositely charged particle states denoted $W^+$ or $W^-$ depending on orientation.

\newpage
The interaction among these fields  can be pictorialized in the following way:
\be W^+ \hskip 3cm \hskip .7cm W^+ \hskip 3.8cm A_\mu \hskip 1.6cm \ee
\vskip -.7cm
\begin{figure}[h]
\centering
\includegraphics{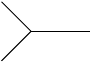}$A_\mu$ \hskip 2cm  \includegraphics{3pt.pdf}$W^+$ \hskip 2cm  \includegraphics{3pt.pdf}$W^-$
\end{figure} 
\vskip -1.7cm
\be
W^- \hskip 3cm  \hskip .7cm A_\mu \hskip 4cm W^-  \hskip 2cm \non
\ee
The three interactions depicted above, $W^+ W^-\rightarrow A_\mu$, $W^+ A_\mu \rightarrow W^+$, and $A_\mu W^-\rightarrow W^-$, all arise from the same interaction vertex and are related to each other by CPT symmetry. 

These interactions can then be manifested as Lie algebra commutators: 
\be [W^+, W^-]=A_\mu \; , \hskip .8cm  \hskip .5cm [W^+, A_\mu]=W^+ \, , \hskip .8cm  \hskip .5cm [A_\mu , W^-]=W^- \; .  \hskip -.5cm\ee
These commutators form precisely the $\mathfrak{su}(2)$ Lie algebra. 

 In addition to the manifestation of interactions as commutators, 
when the singularities are such that the exceptional divisor contains more than  a single cycle, the intersections between the cycles can be interpreted as interactions between fields wrapped or reduced on those cycles \cite{Mayr:1998tx}.

\section{Lie Algebras of the Third Kind (LATKes)} \label{L3}
In this section we define the algebraic objects that are related to codimension 6 orbifold singularities in a way analogous to the relation between Lie algebras and codimension 4 singularities.

Recall from Section \ref{CIsection} that in the correspondence between codimension 4 singularities and Lie algebras, the intersection numbers of pairs of cycles in the exceptional divisor determine the Lie algebra commutators via the CI Relations (equations \ref{CI}). In particular, the intersection numbers enter the following commutators:
\bea \label{CI2}&&[H_i, X_j]=-I_{ij}X_j \; ; 
\\&&  [H_i, Y_j]=I_{ij}Y_j \; .\non
\eea
Also recall the physics interpretation of commutators and intersections as interactions between two fields (Section \ref{GE}). 

Now, for codimension $2n$ singularities for any $n$, the components of the exceptional divisor are $(2n-2)$--cycles, and the intersection of a pair of those has dimension 
\be \dim (C_1 \cap C_2)=\dim C_1+\dim C_2-2n=2n-4.\ee 
Therefore, for codimension 6 orbifolds, the components of the exceptional divisor are 4--cycles, and the intersection of any pair $C_1, C_2$ of 4--cycles does not yield a number but a two-dimensional space:
\be \dim (C_1 \cap C_2)=4+4-6=2.\ee 

However, the intersection of {\it three} 4--cycles in a six dimensional ambient space is zero--dimensional, yielding intersection numbers with three indices: $I_{ijk}$.  On the physics side, these triple intersections should lead to interactions among three fields. 

Bringing together all the above leads us to a generalization of the CI Relations to the codimension 6 case. Equations (\ref{CI2}) become commutators of three objects:
\bea \label{CI3}&&[A_i, B_j, X_k]=-I_{ijk}X_k \; ; 
\\&&  [A_i, B_j, Y_k]=I_{ijk}Y_k \; .\non
\eea
At this stage, we do not yet know what the $A_i$, $B_j$, $X_k$, or $Y_k$ are. However, equation (\ref{CI3}) provides the fundamental ingredient in the algebraic objects we have been searching for: a commutator with three entries. It is now clear that the sought-after algebraic  objects are the following natural generalizations of Lie algebras:

\begin{defin} A {\textit {\textbf {Lie algebra of the third kind (a "LATKe")}}} $\L$ is a vector space 
equipped with a  commutator of the third kind, which is a trilinear anti-symmetric map  
\be [\hskip .1cm\cdot \hskip .1cm , \hskip .1cm\cdot \hskip .1cm, \hskip .1cm\cdot \hskip .1cm]: \Lambda ^3 \L \rightarrow \L \ee
  that satisfies the  Jacobi identity of the third kind:
\be \label{3JI} [X,Y,[Z_1, Z_2, Z_3]] =[[X,Y,Z_1],Z_2, Z_3]+[Z_1, [X,Y,Z_2], Z_3]+[Z_1, Z_2, [X, Y, Z_3]]
\ee
for $X, Y, Z_i\in \L$. 
\end{defin}

We will refer to the commutator of the third kind as a "LATKe commutator", and to the Jacobi identity of the third kind as the "LATKe Jacobi identity." The LATKe Jacobi identity is simply a generalization of the standard Jacobi identity, which says that the adjoint action is a derivation. 

For codimension $2n$ singularities, it is $n$--fold intersections of $(2n-2)$--cycles that give intersection numbers. So the CI Relations for that case have commutators of $n$ objects, leading to the following generalization of a LATKe:

\begin{defin} A {\textbf {\textit  {Lie algebra $\L$ of the $n$-th kind (a "LAnKe")}}}   is a vector space equipped with a  commutator  of the $n$-th kind, which is an n-linear, totally antisymmetric  map 
\be [\cdot, \cdot , \hskip .7cm , \cdot] : \wedge ^n \L \rightarrow \L ~\ee
that satisfies the Jacobi identity of the $n$-th kind:
\be [X_1, \ldots X_{n-1},[Z_1, \ldots Z_n]]=\sum _{i=1}^n [Z_1,\ldots , [X_1, \ldots X_{n-1},Z_i], \ldots Z_n]~,\ee
for $X_i, Z_j \in \L$.
\end{defin}

Before we go any further, we should satisfy ourselves that LATKes  actually exist. Therefore, we now construct an example. 
\vskip .3cm
\no \begin{underline} {\tt The algebra $\L_3$}\end{underline} 
Recall that in the codimension 4 case, each component of the exceptional divisor corresponds to a node in the Dynkin diagram of the corresponding Lie algebra. Therefore, each component corresponds to a simple root of the Lie algebra.
Furthermore, the simplest singularity is $\C^2/\Z_2$, whose blow-up is a single $\P ^1$, and the simplest (non-abelian) Lie algebra is the corresponding $\mathfrak{su}(2)$. 

It is reasonable to presume that similarly,  the simplest example of a LATKe  would correspond to the codimension 6 orbifold singularity with the simplest exceptional divisor. Consider the singularity $\cz$, where the action of a generator $\epsilon$ of $\Z_3$ is given by
\be \label{epsilonaction}\epsilon : (x, y, z)\longmapsto (\epsilon x , \epsilon y, \epsilon z)\; , \mbox{ where }\;  \epsilon ^3=1 , (x,y,z)\in \C^3.\ee
The blow-up of this singularity is a single $\P^2$ (see Appendix \ref{blow}). The cycle $\P ^2$ corresponds to a node in the yet-to-be-defined   Dynkin diagram of $\L_3$. That node, in turn, should correspond to a yet-to-be-defined "root" of $\L _3$. 

We must now define the notion of a root for LATKes. For standard Lie algebras, a root $\alpha$ is in the dual space of the Cartan subalgebra $\h $: 
\be \alpha : \h \longrightarrow \C \; .\ee
So, to define a  root for $\L$, we must first define a  Cartan subalgebra for $\L$. 

In the standard Lie algebra case, one can think of the Cartan subalgebra as a set of operators acting on the Lie algebra $\g$: for a root vector $X_\alpha \in \g$ and for $H\in \h$, we have
\bea \hskip 2cm H&:&X_\alpha \longmapsto [H, X_\alpha ]=\alpha (H)X_\alpha ~.\eea
For a LATKe $\L$, there is no natural action of a subalgebra. However, if $\h _\L \subset \L$ is a  subalgebra of $\L$ (a subspace closed under the commutator), then $\Lambda ^2 \h _\L$ does act on $\L$ naturally:  if $H_1, H_2\in \h_\L$ and $X\in \L$, then $H_1\wedge H_2\in \Lambda ^2\h_\L$ acts on $\L$ via:
\bea \label{CartanAct}
H_1\wedge H_2 &:& X \longmapsto [H_1, H_2, X] \; .
\eea
We may now define $\h_\L \subset \L$ the same way it is defined for Lie algebras: 
\begin{defin} A {\textit {\textbf {Cartan subalgebra}}} $\h _\L$ of $\L$ is a maximal commuting subalgebra of $\L$ such that $\Lambda ^2\h _\L$ acts diagonally on $\L$. 
\end{defin}
Finally, following the standard Lie algebra case, we have 
\begin{defin} Let $\L$ be a LATKe and let $\h _\L$ be a  Cartan subalgebra of $\L$. A {\textit {\textbf {root}}} $\alpha$ of $\L$  is a map in the dual space of $\Lambda ^2 \h _\L$:
\be \alpha : \Lambda ^2 \h _\L \longrightarrow \C \; . \ee
\end{defin}
Now we may also define a generalization for Dynkin diagrams:
\begin{defin} \label{3Dynk} A {\textit {\textbf {Dynkin diagram}}} of $\L$ is a graph with a vertex for each simple positive root of $\L$.
\end{defin}
By "positive," we imply that we have chosen a direction in $(\Lambda ^2\h _\L)^*$ and ordered the roots, as is done for conventional Lie algebras. Note that while this  definition  for a Dynkin diagram may appear to be only a partial  one because it says nothing about edges, it will become clear later that this definition is 
complete.

We now construct a LATKe with a single root and a single node in its  Dynkin diagram corresponding to the cycle $\P ^2$ in the blow-up of $\cz$.
 
For the  root  space to be one-dimensional, the  Cartan subalgebra must be two-dimensional. Let $\h _\L = \{H_1, H_2 \}$ so that $\Lambda ^2 \h _\L$ is spanned by $H_1 \wedge H_2$, and introduce candidates for one positive and one negative root vector, named $X$ and $Y$, respectively. All of the above lead to the following two equations:
\bea \;  [H_1, H_2, X]&=&\alpha (H_1\wedge H_2)X \;  ;  \\
\; [H_1, H_2, Y]&=&-\alpha (H_1\wedge H_2)Y \;  ,
\eea
where $\alpha$ is the single root associated with $\Lambda ^2\L$. Normalizing $H_1$ and/or $H_2$ so that $\alpha (H_1\wedge H_2) =1$ (we should really be normalizing to  $\alpha (H_1\wedge H_2) =-I_{111}$, where $I_{111}$ is the triple intersection number of the $\P ^2$; however, the resulting algebra would be equivalent) gives
\bea \label{rootcomm}\; \hskip -2cm [H_1, H_2, X]&=& X \; ;  \\
  \; \hskip -2cm [H_1, H_2, Y]  &=&  -Y \; . \non
\eea

All that is left is determining $[H_1, X, Y]$ and $[H_2, X, Y]$ such that the LATKe Jacobi identity would be satisfied. This identity requires, among other things, that $[H_i, X, Y], i=1,2$ are zero roots:
\bea [H_1, H_2,[H_i, X, Y]]&=&[[H_1, H_2, H_i], X, Y]+[H_i, [H_1, H_2, X], Y]\non \\
&& \hskip 4cm + [H_i, X, [H_1, H_2, Y]]   \non \\  
&=&0+[H_i, X, Y]+[H_i, X, -Y]=0
\eea
Therefore, $[H_i, X, Y]\in \h _\L$. 

We now restrict the possibilities for $[H_i, X, Y]$ by limiting ourselves to {\sl simple} LATKes, which we now define (recall that in the codimension 4 case, all the Lie algebras corresponding to orbifolds are simple).
 First we need the notion of ideals. 
\begin{defin} An {\textit {\textbf {ideal}}} of $\L$ is a subalgebra ${\cal I}$ that satisfies
\be  [\L, \L, {\cal I}]\subset {\cal I} \; .\ee
\end{defin}
\begin{defin} A  LATKe is {\bf {\textit {simple}}} if it is non-abelian and has no non-trivial ideals. 
\end{defin}

Requiring that our example would be simple means that $[H_1, X, Y]$ and $[H_2,X,Y]$ must be linearly independent, otherwise there would be a non-trivial ideal. It also requires that $[H_i, X, Y]$ and $H_i$ be linearly independent for each $i$. 

So $[H_i, X, Y]$ are linearly independent combinations of $H_1$ and $H_2$. Any such linear combination can be transformed into another by  linear redefinitions that do not affect the two commutators we already fixed in equation (\ref{rootcomm}). Therefore, we may make a choice and we let
\bea \label{restcomm}\; \hskip -2cm [H_1,  X, Y]&=& H_2 \; ;  \\
  \; \hskip -2cm [H_2, X, Y]  &=&  H_1 \; . \non
\eea
No more commutators need to be determined -- all other commutators are defined by the requirement that the commutator map is trilinear and totally antisymmetric. One can easily verify that the LATKe Jacobi identity  is satisfied and that there are no non-trivial ideals. Therefore, we have 
\begin{eg} \label{L3eg} The four-dimensional algebra spanned by $\{ H_1, H_2, X, Y\}$ with the  commutators 
\bea \label{L3com}\; \hskip -2cm [H_1, H_2, X]&=& X \;  ,\non \\
  \; \hskip -2cm [H_1, H_2, Y]  &=&  -Y \;  ,  \\
\; \hskip -2cm [H_1,  X, Y]&=& H_2 \;  , \non \\
  \; \hskip -2cm [H_2, X, Y]  &=&  H_1 \;  ,\non
\eea
is a simple LATKe. We name it $\L _3$. Its Cartan subalgebra is 2-dimensional, 
its root space is 1-dimensional, and its Dynkin diagram consists of a single point:
\[ \bullet  \]
\vskip .5cm
\end{eg}

\section{Uniqueness of the LATKe}\label{1-proof}
We have constructed the LATKe $\L_3$ corresponding to the simplest codimension 6 orbifold singularity, $\cz$. 
If the definition of LATKes has planted seeds for  a generalization of the Du Val-Artin correspondence, then there should be a simple  LATKe  for each of  the orbifolds $\C ^3/\G$, where $\G$ is any discrete, finite subgroup of $SU(3)$ (these $\G$'s are classified in \cite{YauYu}). 

Our goal now, therefore, is to classify all simple LATKes. In doing so, we would also learn more about their structure, which should help us in eventually constructing a physics theory appropriate to these algebras.

We begin  by recalling that we can think of $\Lambda ^2 \h _\L$ as a set of operators acting on $\L$, as in equation (\ref{CartanAct}). More generally, we note that $\Lambda ^2 \L$, not just $\Lambda ^2 \h _\L$, acts on $\L$ as follows: if $X\wedge Y \in \Lambda ^2 \L$ then its action on $Z\in \L$ is given by
\bea 
\label{adaction} X\wedge Y &:& Z \longmapsto [X, Y, Z] \; .
\eea
This action is called {\textit {\textbf {the adjoint action}}} on $\L$ and we denote it ad$(X\wedge Y)$. 

The LATKe Jacobi identity, equation (\ref{3JI}), indicates that $\Lambda ^2 \L$ acts as an {\textit {inner derivation}} of $\L$:
\begin{defin} An operator $D$ on $\L$ is a {\textit {\textbf {derivation}}} of $\L$ if it satisfies
\be D [X,Y,Z]=[DX,Y,Z]+[X,DY,Z]+[X,Y,DZ] \; . \ee
The set of derivations of $\L$ forms a Lie algebra with bracket given by
\be \label{bracketIDer} [D_1, D_2]=D_1D_2-D_2D_1~.\ee
\end{defin}
\begin{defin} An {\textit {\textbf {inner derivation}}} of $\L$ is a derivation of $\L$ which is a linear combination of operators of the form given in equation (\ref{adaction}). 
\end{defin}
These two definitions are analogous to the ones given for conventional Lie algebras. 

We shall denote the algebra of derivations of $\L$ by Der$(\L)$ and the inner derivations by IDer$(\L)$. 
\begin{defin} \label{gL} For any LATKe $\L$, the {\textit {\textbf {Lie algebra of $\L$}}} is the space $IDer(\L)$ with the commutator given by equation (\ref{bracketIDer}); we  denote this Lie algebra $\g _\L$.
\end{defin}

Note that the adjoint action defined before, given by  
\be \label{adgen}\mbox {ad}: X\wedge Y \in \Lambda ^2 \L\longmapsto \mbox{ad}(X\wedge Y)\in IDer(\L) ~,\ee
where
\be \label{adspec}\mbox{ad}(X\wedge Y)(Z)=[X,Y,Z]~,\ee
satisfies
\be \label{intertw} [D, \mbox{ad}(X\wedge Y)] = \mbox{ad}(D(X\wedge Y)) \hskip 1cm\forall D\in IDer(\L),
\ee
and is a surjective morphism of representations of IDer$(\L)$.

It will be convenient to rewrite the  LATKe Jacobi identity in the following form:
\bea \label{JII}[X_1, X_2 , [X_3, X_4, X_5]]-[X_3, X_4,[X_1, X_2, X_5]]=\hskip 2cm \non \\ \hskip 2cm = [[X_1, X_2, X_3], X_4, X_5]+[X_3, [X_1, X_2, X_4], X_5]~,
\eea
so that
\be \label{gLcomm} [X_1 \wedge X_2, X_3 \wedge X_4]=[X_1, X_2, X_3]\wedge X_4 + X_3\wedge [X_1, X_2, X_4]~,
\ee
or equivalently
\be \label{gLcommad}[\mbox {ad}(X_1 \wedge X_2),\mbox {ad}(X_3 \wedge X_4)]=\mbox {ad}([X_1, X_2, X_3]\wedge X_4)+\mbox {ad}(X_3\wedge [X_1, X_2, X_4])~.
\ee

\vskip .3cm

\no \begin{underline}{\tt  $\L$ simple $\Longrightarrow \ggg$ semi-simple}\end{underline} We have shown that every LATKe $\L$ has associated with it a Lie algebra $\ggg$. We now prove a series of lemmas culminating in the result that if $\L$ is a simple LATKe then its Lie algebra $\ggg$ is semi-simple. Consequently,
 we will be able to utilize the well-known classification of semi-simple Lie algebras as a tool for classifying the simple LATKes.

\begin{lemma} \label{Lirrep} If $\L$ is simple then $\L$ is irreducible as a representation of IDer$(\L)$.
\end{lemma}
\no {\bf Proof} Let $W\subset \L$ be an invariant subspace, i.e. IDer$(\L): W\rightarrow W$. Then $[\L, \L, W]\subset W$ (see equations (\ref{adgen}) and (\ref{adspec})). So $W$ is an ideal of $\L$. Since $\L$ is simple, $W=0$ or $\L$. 
\begin{lemma} \label{idealIDer} IDer$(\L)$ is an ideal of Der$(\L)$.
\end{lemma}
\no {\bf Proof} Let  $D\in$ Der$(\L)$ and let $X\wedge Y\in  \Lambda ^2 \L$ so ad$(X\wedge Y)\in$ IDer$(\L)$. Then we have 
\be  [D, \mbox{ad}(X\wedge Y)]\cdot Z = D[X,Y,Z] -[X,Y,DZ]=[DX,Y,Z]+[X,DY,Z]~, \ee
so
\be [D,\mbox{ad}(X\wedge Y)]=\mbox{ad}(DX\wedge Y)+\mbox{ad}(X\wedge DY)~, \ee
which is  an inner derivation.
\begin{fact} The space $\L$ is a representation of Der$(\L)$ and it is faithful by definition.
\end{fact}
\begin{lemma} If $\L$ is simple then $\L$ is irreducible as a representation of Der$(\L)$.
\end{lemma}
\no {\bf Proof} Any subspace $W\subset \L$ that is invariant under Der$(\L)$ is also invariant under its ideal IDer$(\L)$. By the proof of Lemma \ref{1-proof}.1, $W=0$ or $\L$. 
\begin{lemma} If $\L$ is simple then Der$(\L)$ is reductive.
\end{lemma}
\no {\bf Proof} This follows from the fact that Der$(\L)$ has  a finite dimensional, faithful, irreducible representation, namely $\L$.
\begin{lemma} \label{derIder}If $\L$ is simple then any derivation of $\L$ is an inner derivation, i.e. Der$(\L)=$IDer$(\L)$.
\end{lemma}
\no {\bf Proof} Let $D$ be any derivation of $\L$. Since Der$(\L)$ is reductive, it has the form of a direct sum of commuting ideals. We already know that IDer$(\L)$ is an ideal of Der$(\L)$. Therefore, all we need to show is that if  $D$ commutes with IDer$(\L)$, then $D=0$.
Assume 
\be [D, \mbox{ad}(X\wedge Y)]\cdot Z=0 \hskip 1cm \forall X,Y,Z\in \L~.\ee
Expand this equation to  
\be D(\mbox{ad}(X\wedge Y)\cdot Z) - \mbox{ad}(X\wedge Y)(D\cdot Z)=D[X,Y,Z]-[X,Y,DZ]=0.\ee
Repeating this for permutations of $X,Y$, and $Z$ and using the definition of derivations, we find that 
\be [DX,Y,Z]=[X,DY,Z]=[X,Y,DZ]=0 \hskip 1.5cm \forall X,Y,Z\in \L~.\ee
Therefore, $DX$ is in the center of $\L$ for any $X\in \L$. But since $\L$ is simple, it has no center. Therefore, $DX=0 ~\forall X$ so $D=0$.  
\begin{lemma} If $\L$ is simple, then the center of IDer$(\L)$ is trivial.
\end{lemma}
\no {\bf Proof} The argument in the proof of Lemma \ref{1-proof}.6 shows that any derivation that commutes with IDer$(\L)$ is zero.
\begin{lemma} If $\L$ is simple then IDer$(\L)$ is semi-simple.
\end{lemma}
\no {\bf Proof}  Since IDer$(\L)$ is an ideal in a reductive Lie algebra, it is itself reductive. Any reductive Lie algebra is the direct sum of a semi-simple part and its center. Since IDer$(\L)$ has no center, it is semi-simple.

\vskip .3cm
So we have shown that IDer$(\L) = \ggg$ is a semi-simple Lie algebra when $\L$ is simple. We show next that there are very strong constraints on the roots and weights of $\ggg$ that substantially limit the number of potential candidates for $\ggg$.  

\vskip .3cm

\no \underline{{\tt  Constraints on roots and weights of $\ggg$}}
 Since $\L$ is a representation of $\ggg$, so is $\Lambda ^2\L$. 
Recall from equations (\ref{adgen}), (\ref{adspec}), and  (\ref{intertw}) that we have a surjective morphism of representations, 
\be \mbox{ad}: \Lambda ^2\L  \longrightarrow \ggg ~.\ee
Therefore, there is a relation between the weights of  $\Lambda ^2\L $ and the roots of $\ggg$, which we shall now explore. 

Let $\h$ be the Cartan subalgebra of $\ggg$, let $H\in \h$, and let $X_{\beta  _i}\in \L$ be the weight vectors of $\L$ with $\beta _i$ denoting weights of $\L$ so that 
\be H(X_{\beta  _i})= \beta  _i(H)   X_{\beta  _i}~.\ee     
Then 
\bea \label{sumweights}[H ,\mbox{ad}(X_{\beta  _i}\wedge X_{\beta  _j})]&=& \mbox{ad}(H( X_{\beta  _i}\wedge X_{\beta  _j})  )  \\
&=&(\beta _i (H)+\beta _j (H))\mbox{ad}(X_{\beta  _i}\wedge X_{\beta  _j})\;  . \nonumber
\eea   
It follows that any root of $\ggg$ has the form $\beta _i + \beta _j$.
However, since ad may have a non-trivial kernel, not all weights $\beta _i + \beta _j$ are necessarily roots. 

Let the highest weight of $\L$ as an irreducible representation of $\ggg$ be denoted $\Lambda$. Now we prove a series of lemmas culminating in the result that $2\Lambda -\alpha=\Lambda + (\Lambda -\alpha)$ is a highest root of $\ggg$, where $\alpha$ is a simple positive root of $\ggg$. We do so by showing that there is a lowering operator $E_{-\alpha}\in \g _\L$ such that if $v_\Lambda$ is the highest weight vector of $\L$, then $v_\Lambda \wedge (E_{-\alpha}v_\Lambda)$ is not in the kernel of ad. 

The first step is to construct a (different) element in $\Lambda ^2 \L$ which is not in the kernel of ad. Let $w_0$ be the element  of the Weyl group of IDer$(\L)$ that takes every positive root to a negative one and vice versa; then $w_0^2=1$.

\begin{lemma} We have ad($v_\Lambda \wedge v_{w_0\Lambda}$)$\neq 0$ and ad($v_\Lambda \wedge v_{w_0\Lambda}$)$\in \h$.
\end{lemma}
\no {\bf Proof} It is easy to see that $v_\Lambda \wedge v_{w_0\Lambda}$ generates the entire space $\Lambda ^2\L$. Therefore, its image under ad must be nonzero, otherwise ad itself would be identically zero, which would contradict that ad is onto IDer$(\L)$.

Now, $\Lambda + w_0\Lambda$ is a root (see equation (\ref{sumweights})). Since $w_0(\Lambda+ w_0\Lambda)=w_0\Lambda + w_0^2\Lambda = w_o\Lambda +\Lambda$ and $w_0$ takes positive roots to negative ones and vice versa, it follows that $\Lambda + w_0\Lambda=0$, or $w_0\Lambda = -\Lambda$. Therefore, ad($v_\Lambda \wedge v_{w_0\Lambda}$)$\in \h$.

\begin{lemma} In each simple ideal of IDer$(\L)$, there is a simple positive root $\alpha$ such that ad($v_\Lambda \wedge (E_{-\alpha}v_{\Lambda})$)$\neq 0$.
\end{lemma}

\no {\bf Proof} We have
\bea \label{notkernel}\mbox{ad}(v_\Lambda \wedge (E_{-\alpha}v_{\Lambda}))\cdot v_{w_0\Lambda}&=& [v_\Lambda , (E_{-\alpha}v_{\Lambda}), v_{w_0\Lambda}] \non
\\  &=& \mbox{ad} ((E_{-\alpha}v_{\Lambda}) \wedge v_{w_0\Lambda})\cdot v_\Lambda \non
\\ &=&\mbox{ad}(E_{-\alpha}(v_{\Lambda}\wedge v_{w_0\Lambda}))\cdot v_\Lambda 
\\ &=& [E_{-\alpha}, \mbox{ad}(v_{\Lambda}\wedge v_{w_0\Lambda})]\cdot v_\Lambda \non
\\ &=& \alpha (\mbox{ad}(v_{\Lambda}\wedge v_{w_o\Lambda}))E_{-\alpha} \cdot v_\Lambda ~.\non
\eea
If we now prove that for each simple ideal of IDer$(\L)$, there is a simple root $\alpha$ such that $\alpha (\mbox{ad}(v_{\Lambda}\wedge v_{w_0\Lambda}))\neq 0$ and $E_{-\alpha} \cdot v_\Lambda \neq 0$, we will have proven the Lemma.  

First, if we have $\alpha (\mbox{ad}(v_{\Lambda}\wedge v_{w_0\Lambda}))\neq 0$, then the identity
\bea \mbox{ad} ((E_{-\alpha}v_{\Lambda}) \wedge v_{w_0\Lambda}) = \alpha (\mbox{ad}(v_{\Lambda}\wedge v_{w_0\Lambda}))E_{-\alpha}~, \eea
which can be deduced from equation (\ref{notkernel}), implies that $E_{-\alpha}v_{\Lambda}\neq 0$. So we need to show only that there is an $\alpha$ such that $\alpha (\mbox{ad}(v_{\Lambda}\wedge v_{w_0\Lambda}))\neq 0$.

Suppose first that IDer$(\L)$ is simple. Then the root space is dual to $\h$ so such an $\alpha$ automatically exists. Now suppose that IDer$(\L)$=$\g_1\oplus \cdots \oplus \g_k$, where the $\g_i$ are ideals in IDer$(\L)$. Then the Cartan subalgebra also  has the form $\h = \h_1 \oplus \cdots \oplus \h_k$, and the root space decomposes to $\Delta = \Delta _1 \oplus \cdots \oplus \Delta_k$. Now, suppose that $\alpha (\mbox{ad}(v_{\Lambda}\wedge v_{w_0\Lambda}))= 0$ for all $\alpha \in \Delta _1$. Then $\mbox{ad}(v_{\Lambda}\wedge v_{w_0\Lambda})\in \h_2\oplus \cdots \oplus \h_k$. Since $v_{\Lambda}\wedge v_{w_0\Lambda}$ generates all of $\Lambda ^2\L$, this means that the image of ad does not contain $\h_1$, which contradicts the surjectivity of ad. 

Therefore, for each $i=1, \ldots, k$, there is a simple root $\alpha$ that satisfies both $\alpha (\mbox{ad}(v_{\Lambda}\wedge v_{w_0\Lambda}))\neq 0$, and $E_{-\alpha} \cdot v_\Lambda \neq 0$, so by equation (\ref{notkernel}), $v_\Lambda \wedge (E_{-\alpha}v_{\Lambda})$ is not in the kernel of ad. \hfill $\Box$

\vskip .3cm
Note that no root of IDer$(\L)$ is higher than $2\Lambda -\alpha$ because there is no available pre-image under ad for such a root in $\Lambda ^2\L$, and ad is surjective. Putting everything together,  we have

\begin{lemma} The highest root $\theta$ of any simple component of $\ggg$ is the sum of $\Lambda$, the highest weight of $\L$, and a next-to-highest weight $\Lambda -\alpha$, where $\alpha$ is a simple positive root of a simple component of $\ggg$:
\be \theta = \Lambda + (\Lambda -\alpha) \ee
or
\be \label{condition}\theta + \alpha = 2\Lambda ~.\ee
\end{lemma}
\vskip .2cm
Only semi-simple Lie algebras which have a faithful irreducible representation $V$ whose highest weight $\Lambda$ satisfies
condition  (\ref{condition}) are potential candidates to be Lie algebras of some $\L$. Our next step, therefore, is to find all semi-simple $\g$ and highest weights $\Lambda$ satisfying this condition. 

As it turns out, this same condition appeared in an entirely different context in Kac's work on Lie superalgebras \cite{kac}, where all pairs of semi-simple Lie algebras $\g$ and irreducible faithful representations $V$ with highest weight $\Lambda$ satisfying this condition are classified. 
The resulting list of pairs is the following.

If $\g$ is not simple then $\g=\mathfrak{so}_4 \simeq \mathfrak{sl}_2 \times \mathfrak{sl}_2   $ and $V$ is the standard four-dimensional (vector) representation.

If $\g$ is simple then the following table constitutes the complete list of all pairs of $\g$ and $V$ that satisfy the condition:
\be \label{list} \begin{array} [h] {llll} \g&\dim V& \dim \g & \dim \Lambda ^2 V \\ \hline
A_1 &3 &3 &3 \\
G_2&7&14&21\\
A_3&6&15&15\\
B_3&8&21&28\\
B_{r\geq 2}\hskip .5cm &2r+1\hskip .3cm &r(2r+1)\hskip .3cm &r(2r+1)\\
D_{r\geq 4}&2r&r(2r-1)&r(2r-1)
\end{array}
\ee
The dimensions in the table uniquely identify each representation, except that for $D_4$ there are three different representations of dimension 8, which are related to each other by the triality symmetry of $\mathfrak{so}_8$.

\vskip .3cm
\no \begin{underline}{\bf {\tt {Further conditions on $\g _\L$ and the construction of all $\L$ }}}\end{underline} 
The condition $2\Lambda = \theta + \alpha$ is necessary, but not sufficient, to ensure that $\g$ is the Lie algebra of some $\L$. 
There is a further requirement:
the intertwining map $\omega : \Lambda ^2 V \rightarrow \g $ must yield a  LATKe commutator in the following way. 
Let $v_1, v_2, v_3\in V$; then we define
\be  [v_1, v_2, v_3]=(\omega (v_1\wedge v_2))\cdot v_3~.\ee
The expression on the right hand side must be antisymmetric in all three variables for it to define a LATKe commutator. Since it is already automatically antisymmetric in the first two variables, the only remaining requirement is 
\be \label{antisym123} (\omega (v_1\wedge v_2))\cdot v_3=-(\omega (v_3\wedge v_2))\cdot v_1 \hskip 1cm \forall  v_1, v_2, v_3\in V \, ,\ee
which is equivalent to requiring
 \be \label{antisym}\omega (v_1\wedge v_2))\cdot v_1=0~ \hskip .5cm \forall v_1, v_2\in V~.\ee

We first prove the following:
\begin{lemma} None of the pairs of $\g$ and $V$ in  table \ref{list} yields a LATKe  commutator. 
\end{lemma}
\no {\bf Proof} There are two steps to the proof. First, we construct explicitly the intertwining map $\omega : \Lambda ^2 V \rightarrow \g $ for each pair of $\g$ and $V$ in table \ref{list}. In each case, the adjoint representation, which is irreducible since $\g$ is simple, appears exactly once in the decomposition of $\Lambda ^2V$. Therefore, by Schur's lemma there is exactly one map $\omega$. Second, we show that $\omega$  does not satisfy equation (\ref{antisym123}) or (\ref{antisym}), so it does not result in a  LATKe commutator. 

We begin with those pairs in table \ref{list} that satisfy $\dim \g = \dim \Lambda ^2V$. We can construct the map $\omega$ for all such pairs simultaneously because all of them have the form $\mathfrak{so}_n$ with $V$ the standard $n$-dimensional representation. That is obvious for $B_r$ and $D_r$, which stand for $\mathfrak{so}_{2r+1}$ and $\mathfrak{so}_{2r}$, respectively. For $A_1$ and $A_3$, recall that $A_1$ stands for $\mathfrak{sl}_2 \simeq \mathfrak{so}_3$ and $A_3$ stands for 
$\mathfrak{sl}_4 \simeq \mathfrak{so}_6$. 

Let $\mathfrak{so}_n$ be spanned by antisymmetric matrices $e_{ab}$, $a,b=1, \ldots , n$, $a<b$, such that
\be \label{eab}\{ e_{ab} \}_{\alpha \beta}=\delta_{a\alpha}\delta _{b\beta}-\delta_{a\beta}\delta _{b\alpha}~, \ee
i.e. $e_{ab}$ has $+1$ in the $ab$-th entry and $-1$ in the $ba$-th entry, with all other entries equal zero; if $a>b$ then $e_{ab}$ is defined by $e_{ab}=e_{ba}$. Let the standard representation $V$ be spanned by the standard basis $\{ e_a \}$, $a=1, \ldots ,n$, where 
\be \label{ea}( e_a )_\alpha=\delta _{a\alpha} ~,\ee
i.e. $e_a$ has $+1$ in the $a$-th entry and zero elsewhere. 
It is straightforward to check that
\bea \label{socomm}[e_{ab} , e_{cd}]&=& \delta _{bc}e_{ad}+ \delta _{ad}e_{bc}-\delta _{ac}e_{bd}-\delta _{bd}e_{ac}~; \\
e_{ab}\cdot ( e_c )&=&\delta _{bc}(e_a)-\delta _{ac}(e_b)~.
\eea

We want to construct $\omega : \Lambda ^2 V \rightarrow \g $ which is intertwining, i.e.
\bea [X, \omega (v_1\wedge v_2)]&=&\omega (X(v_1\wedge v_2)) \non \\ 
&=&\omega ((X\cdot v_1)\wedge v_2+v_2\wedge (X\cdot v_2)) \hskip .3cm \forall X\in \g, v_i\in V  ,
\eea
so we require
\be \label{twining} [e_{ab}, \omega (e_c\wedge e_d)]=\omega ((e_{ab}\cdot e_c)\wedge e_d+e_c\wedge (e_{ab}\cdot e_d))\hskip .3cm \forall a,b,c,d ~.\ee

We now show that the map defined by 
\be \label{thepsi}\omega (e_a\wedge e_b)=e_{ab} \ee
satisfies this property. 

Computing the right hand side of equation (\ref{twining}) for this map gives
\bea  \omega ((e_{ab}\cdot e_c)\wedge e_d+e_c\wedge (e_{ab}\cdot e_d))= \hskip 5cm \non \\
  \hskip 2cm =\omega ((\delta _{bc} (e_a)-\delta _{ac}(e_b))\wedge e_d)+\omega (e_c \wedge (\delta _{bd}(e_a)-\delta _{ad}(e_b)))\\
 = \delta_{bc}e_{ad}-\delta_{ac}e_{bd}+\delta_{bd}e_{ca}-\delta_{ad}e_{cb}~. \hskip 3.7cm \non
\eea
Comparing this with equation (\ref{socomm}) proves that equation (\ref{thepsi}) is the desired map. 

Now we check whether the condition (\ref{antisym}) is satisfied. It is not: 
\be \omega (e_a\wedge e_b)\cdot e_a=e_{ab}\cdot e_a=-e_b\neq 0\; .\ee
Therefore, no LATKe commutator arises from these $\mathfrak{so}_n$'s. 

Now we construct the map $\omega$ for the pair $\g _2$ with $\dim V=7$. As is well-known, $\g_2$ can be realized as the Lie algebra of derivations of the octonions, and the standard 7-dimensional representation of $\g_2$ is given by its action on the (pure imaginary) octonions \cite{Baez:2001dm, FultonHarris}. 

Let $\{ e_i \}$, $i=1, \ldots , 7$ be a basis for the pure imaginary octonions, and for our representation $V$. Define the map $\omega : \Lambda ^2 V\rightarrow \g_2$ by
\be \label{OctonDeriv}\omega (e_i\wedge e_j)\cdot e_k = [[e_i, e_j],e_k]-3((e_ie_j)e_k-e_i(e_je_k))~, \ee
where $[x,y]=xy-yx$, with the multiplication being the one in the octonion algebra, given for example by the standard Fano plane \cite{Baez:2001dm, FultonHarris}. (Recall that the octonions are not associative, so the second term above does not in general vanish; also, this second term is antisymmetric under $e_i\leftrightarrow e_j$, as is the first term). It can be shown \cite{Schafer} that every $\omega (e_i\wedge e_j)$ is a derivation of the octonions, and that every derivation of the octonions is a linear combination of derivations of the form $\omega(e_i\wedge e_j)$. Therefore, the map $\omega$ of Equation (\ref{OctonDeriv}) is onto $\g_2$. Also, it is straightforward to see that $\omega$ is intertwining. Since $\g_2$ appears exactly once in the decomposition $\Lambda ^2 V=\g _2\oplus V$, the map $\omega$ is unique. 

As before, we now check whether the condition  (\ref{antisym}) is satisfied. It is enough to show it is not satisfied in one case:
\be \omega (e_1\wedge e_2)\cdot e_1 = 4e_2\neq 0~.\ee

Therefore, no LATKe commutator arises from $\g _2$.

We have now but one more pair to check:  $B_3=\mathfrak{so}_7$ with $\dim V=8$, where $V$ is the spin representation. The representation $\Lambda ^2 V$ decomposes into $\mathfrak{so}_7\oplus V_7$, where $\mathfrak{so}_7$ stands for the adjoint representation and $V_7$ is the standard seven-dimensional representation of $\mathfrak{so}_7$. Since $\mathfrak{so}_7$ appears once in the decomposition, there is a unique intertwining map $\omega : \Lambda ^2V \rightarrow \mathfrak{so}_7$. We now explicitly construct this map and show that it does not yield a LATKe commutator. 

We first construct the spin representation of $\mathfrak{so}_7$ explicitly, following \cite{FultonHarris}. Let $C(V_7, Q)$ be the Clifford algebra generated by $V_7$ with the quadratic form 
 
 \be Q= \left ( \begin{array} [h] {ccc} 0&I_3 &0\\ I_3&0&0\\ 0&0&1 \end{array}\right )~,
 \ee
 where $I_3$ is the $3\times 3$ identity matrix. Then $\mathfrak{so}_7$ is a Lie subalgebra of $C(V_7, Q)$  via the embedding 
 \be \label{embedding} \psi \cdot \phi ^{-1}:  \mathfrak{so}_7\rightarrow C(V_7, Q)~,
 \ee
 where $\phi: \Lambda ^2V_7\rightarrow \mathfrak{so}_7$ is an isomorphism given by
 \be \phi _{a\wedge b}(v)=2(Q(b,v)a-Q(a,v)b), \hskip 1cm v\in V_7~,
 \ee
 and $\psi: \Lambda ^2V_7 \rightarrow C(V_7, Q)$ is an embedding given by
 \be \psi : a\wedge b \mapsto a\cdot b -Q(a,b)~.\ee
 
 Decompose $V_7$ into $W\oplus W'\oplus U$, where $W$ and $W'$ are three--dimensional isotropic subspaces and $U$ is a one dimensional subspace orthogonal to them. Then there is an action of $C(V_7, Q)$ on $\Lambda ^\bullet W = \sum_{i=0}^3 \Lambda ^iW$, whose restriction to $\mathfrak{so}_7$ will be the spin representation. It is given as follows.
  
  Let $\zeta \in \Lambda ^\bullet W$, and let $l:V_7\rightarrow \mbox{End}(\Lambda ^\bullet W )$ be given by:
  \bea l(w)\cdot \zeta &=& w\wedge \zeta \hskip 1.5cm w\in W~, \non \\
  l(w')\cdot \zeta &=& D_{w'} (\zeta) , \hskip 1cm w'\in W'~, \\
  l(e_0)\cdot \zeta &=& \left \{ 
\begin{array}{rr}
\zeta \hskip 1cm \zeta \in \Lambda ^{even}W \\
-\zeta \hskip .8cm \zeta \in \Lambda ^{odd}W 
\end{array} \right. , \hskip .6cm e_0\in U, \; Q(e_0, e_0)=1. \non 
  \eea 
  where
  \bea &&D_{w'} (1)=0, \nonumber \\
  &&D_{w'}(w)=2Q(w,w')  , \\
  &&D_{w'}(w_1\wedge w_2)
   = -\phi_{w_1\wedge w_2}(w') , \nonumber \\
  &&D_{w'} (w_1\wedge w_2\wedge w_3)=D_{w'}(w_1\wedge w_2)\wedge w_3+ D_{w'}(w_3)w_1\wedge w_2  . \nonumber
  \eea

  The map $l$ gives the action of $V_7$ on $\Lambda ^\bullet W$ from which the action of $C(V_7, Q)$ on $\Lambda ^\bullet W$ can be deduced. This action restricted to $\mathfrak{so}_7$ is the spin representation. 
  
  The action of $\mathfrak{so}_7$ on $\Lambda ^2(\Lambda ^\bullet W)$ is then given by
   \be X(\zeta _1\wedge \zeta _2)=(l(X)\cdot \zeta _1)\wedge \zeta _2 + \zeta _1 \wedge (l(X)\cdot \zeta _2) ~,  \hskip 1cm X\in \mathfrak{so}_7 ~. \ee

  Now we turn to constructing the map $\omega$ satisfying the intertwining condition, i.e. 
  \be [X, \omega (\eta)]=\omega (X(\eta)) \hskip 1cm \forall X\in \mathfrak{so}_7, \eta \in \Lambda ^2(\Lambda ^\bullet W).
  \ee
  If $\eta \in \Lambda ^2(\Lambda ^\bullet W)$ is a unique (up to scalar) element satisfying 
  \be H(\eta)=\alpha (H)\eta \hskip 1cm \forall H\in \mathfrak{h}\subset \mathfrak{so}_7 ~,
  \ee
  where $\mathfrak{h}$ is the Cartan subalgebra of $\mathfrak{so}_7$ and $\alpha \in \mathfrak{h}^*$ is a root, then since $\omega$ is onto $\mathfrak{so}_7$ and is intertwining, $\omega (\eta)$ must be a non-zero multiple of the $\alpha$ root vector in $\mathfrak{so}_7$.

   Let $\{ e_i\}$, $i=1,2,3$ denote a basis for $W$, let $\{ e_{i+3}\}$, $i=1,2,3$ denote a basis for $W'$, and let $e_0$ be a basis for $U$. With the notation given in Appendix \ref{omegaB3}, let $H_1=\mu_{11}$, $H_2=\mu_{22}$, and $H_3=\mu_{33}$ form the Cartan subalgebra.  
  Consider the element $1\wedge e_1\in \Lambda ^2(\Lambda ^\bullet W)$. Then
  \be H_1(1\wedge e_1)=0~; \hskip 1cm H_2(1\wedge e_1)=-1\wedge e_1~; \hskip 1cm H_3(1\wedge e_1)=-1\wedge e_1~.\ee
  One can check that $1\wedge e_1$ is the only element in $\Lambda ^2(\Lambda ^\bullet W)$ (up to scalar) with these eigenvalues. In  $\mathfrak{so}_7$, the unique element (up to scalar) with these eigenvalues is $\rho _{23}$ (see Appendix \ref{omegaB3}). Therefore, we may set 
  \be \label{startomega}\omega (1\wedge e_1)=\rho_{23}~.\ee
  
The rest of $\omega$ is fully  determined by repeated applications of the intertwining condition and the action of $\mathfrak{so}_7$ on $\Lambda ^2(\Lambda ^\bullet W)$. We provide the resulting $\omega$ in Appendix \ref{omegaB3}. 
  
  It is now straightforward to check that the map $\omega$ does not satisfy Equation (\ref{antisym123}):
  \be \omega(f_1\wedge f_5)\cdot f_4=-2f_1 \neq -\omega(f_1\wedge f_4)\cdot f_5 = f_1~,\ee 
  where $\{ f_k \}, k=1, \ldots, 8$ is the basis for $\Lambda ^\bullet W$ given in Appendix \ref{omegaB3}. 

Therefore, we do not get a LATKe from any of the simple Lie algebras in table  \ref{list}, and the proof of the lemma is complete. \hfill $\Box$

\vskip .3cm
We are now left with only one candidate: $\mathfrak{so}_4$ with its standard representation. Recall that $\mathfrak{so}_4$ is not a simple Lie algebra but has two simple factors, $\mathfrak{so}_4\simeq \mathfrak{sl}_2 \times \mathfrak{sl}_2$. So its adjoint representation is not irreducible, Schur's lemma does not apply, and the intertwining map $\omega$ constructed above for  $\mathfrak{so}_n$ is not unique for this case. We can construct another one. 

We will show that the other map  does lead to a  LATKe commutator and in fact yields the LATKe $\L_3$ which we constructed in Section \ref{L3}. 

Using the same notation as before, the basis for $\mathfrak{so}_4$ is
\be \{ e_{12}, e_{13},e_{14},e_{23},e_{24},e_{34} \} \ee
and the basis for $V$ is 
\be \{ e_1, e_2, e_3, e_4 \}~.\ee
We define $\phi: \Lambda ^2 V \rightarrow \g $ explicitly by
\bea \phi (e_1\wedge e_2)&=&e_{34}, \non \\ 
\phi (e_1\wedge e_3)&=&-e_{24}, \non \\ 
\phi (e_1\wedge e_4)&=&e_{23}, \non \\ 
\phi (e_2\wedge e_3)&=&e_{14}, \label{phi4}\\ 
\phi (e_2\wedge e_4)&=&-e_{13}, \non \\ 
\phi (e_3\wedge e_4)&=&e_{12}. \non 
 \eea
 It is straightforward to check that this map satisfies the intertwining condition, equation (\ref{twining}). 
 
 We now use the map $\phi$ to construct the  LATKe commutator. The dimension of the candidate for $\L$ is $\dim V=4$, so there are only four commutators to calculate:
 \bea \; [e_1, e_2, e_3]&=&\phi(e_1\wedge e_2)\cdot e_3=e_{34}\cdot e_3=-e_4 \; ; \non \\
 \; [e_1, e_2, e_4]&=&\phi(e_1\wedge e_2)\cdot e_4=e_{34}\cdot e_4=e_3 \; ; \non \\
\;  [e_1, e_3, e_4]&=&\phi(e_1\wedge e_3)\cdot e_4=-e_{24}\cdot e_4=-e_2 \; ;  \label{L3pretty}\\
 \; [e_2, e_3, e_4]&=&\phi(e_2\wedge e_3)\cdot e_4=e_{14}\cdot e_4=e_1 \; .\non 
 \eea
 With the following change of variables
 \be X=\frac{1}{\sqrt 2}(e_1+ie_4)\; ; Y=\frac{1}{\sqrt 2}(e_1-ie_4)\; ; H_1=-ie_2\; ; H_2=e_3 \; ,
 \ee
 we see that this algebra is precisely $\L_3$ as given in equations (\ref{L3com}).
 
 This construction of $\L_3$ generalizes to a LAnKe $\L _n$, as shown in Appendix \ref{lanke}.
 
 We summarize the results of this section in one theorem:
 \begin{theorem} \label{bigtheorem} There is precisely one simple LATKe, namely $\L_3$; it is four dimensional, it corresponds to the singularity $\cz$, its Lie algebra $\g_{\L_3}$ is $\mathfrak{so}_4\simeq \mathfrak{sl}_2\times \mathfrak{sl}_2$, and its Dynkin diagram consists of one node. Its commutators are given by equation (\ref{L3com}) or (\ref{L3pretty}).
 \end{theorem}
  
 We comment that the theorem is equivalent to the following statement. Let $V$ be a vector space with a non-degenerate symmetric  bilinear form so that $\Lambda ^2V\simeq \mathfrak{so}(V)$. There is  a natural action of $\Lambda ^2V$ on $V$ which induces a map $\omega _3: V\otimes V  \otimes V \rightarrow V$, antisymmetric in the first two factors. Let $\bar \omega _3: \Lambda  ^3V\rightarrow V$ be the fully-antisymmetric version of $\omega _3$. The theorem says that, if we assume that the action of $\mathfrak{so}(V)$ on $V$ was irreducible, then there is a unique vector space $V$ for which $\bar \omega _3$ is both non-trivial and satisfies the Jacobi identity of the third kind (equation (\ref{3JI})). That vector space is four dimensional and $\Lambda ^2V$ consists of skew-symmetric 4-matrices.

\section{The Physics of LATKes}\label{physics}
We have now reached an important juncture. Having found the unique, simple LATKe, we ask ourselves: are there any applications of the LATKe to physics? 

The first type of physical theory that comes to mind in applying Lie algebras to physics is gauge theory. Can we generalize gauge theory for LATKes?

The answer is "yes," and we do so in the context of particle physics.

\subsection{LATKe representations and LATKe groups}

\vskip .3cm 
\no \begin{underline}{\texttt  {Representations}}\end{underline} Whenever Lie algebras are applied to particle physics, particles are viewed as basis vectors of representations of the Lie algebra. For example,  gauge fields form the adjoint representation of the gauge group; quarks form the three-dimensional representation of color-$SU(3)$; electrons form the two-dimensional spin representation of $SU(2)$ \cite{Georgi:1982jb}. Therefore, in order to apply LATKes to particle physics,  we must define a "representation" for LATKes. 

To do so, we review the standard Lie algebra case. A representation of a Lie algebra $\g$ is a map from $\g$ to operators on some vector space $V$: 
\be \rho:\g \longrightarrow \mbox{End}(V) ~,\ee
 and it respects the commutator in the following way:
\be \label{repcond2}[\rho (X), \rho (Y)]=\rho ([X, Y]) ~.\ee
A particular representation that utilizes the commutator in a natural way is
the adjoint representation, given by the following map:
\bea \mbox{ad}&:&\g \longrightarrow \mbox{End}(\g) \\
\mbox{ad}(X)&:&Y\longmapsto [X,Y] ~.
\eea
This map satisfies the condition
\be \label{repsrel}[\mbox{ad}(X), \mbox{ad}(Y)] =\mbox{ad}([X,Y])\; , \ee
which is equivalent to the  standard Jacobi identity. The condition in equation (\ref{repcond2}) is a generalization of the relation given in equation (\ref{repsrel}) and it reduces to it when $\rho = $ad. 

Now, we define a representation for a LATKe. We begin by defining the analog of the adjoint representation: it is also a map that utilizes the commutator in a natural way, and we have in fact seen it before (equation (\ref{adaction})):
\bea \label{adL}\mbox{ad}&:&\L\wedge \L\longrightarrow \mbox{End}(\L)\\
\mbox{ad}(X\wedge Y) &:&Z\longmapsto [X,Y,Z]~.
\eea
The map ad satisfies  the condition
\be \label{3reps}[\mbox{ad} (X_1 \wedge X_2),\mbox{ad}( X_3 \wedge X_4)]=\mbox{ad}([X_1, X_2, X_3]\wedge X_4) +\mbox{ad}( X_3\wedge [X_1, X_2, X_4])~,
\ee
which is equivalent to the  LATKe Jacobi identity. 

If we generalize equations (\ref{adL})  and (\ref{3reps}), we have
\begin{defin} A {\textit {\textbf {representation}}} of a LATKe $\L$ is a map
\be \label{latkerep} \rho: \Lambda ^2 \L\longrightarrow \mbox{End}(V) \ee
for some vector space $V$ subject to the condition
\be  [\rho (X_1 \wedge X_2),\rho( X_3 \wedge X_4)]=\rho([X_1, X_2, X_3]\wedge X_4) +\rho( X_3\wedge [X_1, X_2, X_4])~.
\ee
\end{defin}
This condition is analogous to equation (\ref{repcond2}) and it generalizes the  LATKe Jacobi identity (see equations (\ref{3JI}), (\ref{JII}) and (\ref{gLcomm})).
\vskip .3cm
\no \begin{underline}{\tt {Groups}}\end{underline}  Another fundamental ingredient whenever Lie algebras are applied to particle physics is the Lie group, which plays the role of a symmetry of the physical system.

Therefore, in order to apply LATKes to particle physics, we also should define a "Lie group of the third kind," or a  "LATKe group," which would be related to the LATKe in a way analogous to the relation between an ordinary Lie group and its Lie algebra. The LATKe group could then play the role of some kind of generalized symmetry in the yet-to-be constructed physics of LATKes. 

Here, we run into trouble: we have found it impossible to generalize the notion of a Lie group to something we might have called a  LATKe group. While there does happen to be a Lie group associated with the LATKe, namely the exponential of $\g _\L$, it is not in any way an exponential of the LATKe itself. So it is not a "LATKe group."

Since we have no LATKe analog of a Lie group, it would be impossible to generalize in a natural way any application of Lie algebras to physics in which the Lie group is an indispensable ingredient. We are limited to applications in which the only necessary mathematical ingredients are those for which we do have a LATKe analog. 

Since we wish to generalize gauge theory for LATKes, we must check whether 
the Lie group itself, which plays the role of the gauge group, is an indispensable ingredient in the construction of  gauge theory. If it is, we would be unable to generalize it for LATKes.
 In the next section we demonstrate that the gauge group {\sl is not} indispensable in gauge theory by rephrasing the original theory of Yang and Mills \cite{Yang:1954ek} so that all group tranformations are re-written as  Lie algebra actions. That sets the stage for a natural generalization of Yang-Mills theory to LATKes, which we construct in Section \ref{ttYM}.

\subsection{Traditional Yang-Mills theory}\label{YMsans}

Let $\psi$ be a wave function describing a field in some representation $\rho$ of a Lie algebra $\g$. Let $\g$ be spanned by basis elements $T^i$ with $i=1, \ldots , \dim \g$ (these are analogues of Pauli spin matrices in the $\mathfrak{su}(2)$ case). A gauge transformation acts via
\be \label{deltapsi}\delta \psi = -2i\Theta ^iT^i\psi \; ,\ee
where each $\Theta ^i$, $i=1, \ldots , \dim \g$ is a space-time dependent field, and each $T^i$ acts on $\psi$ via the representation $\rho$. 

Invariance under such transformations is preserved only if we require derivatives of $\psi$ to appear in the combination
\be \label{covder} (\partial _\mu -ig B_\mu)\psi \; ,\ee
where $g$ is a coupling constant and 
\be B_\mu = 2b_\mu^j T^j \ee
with $b_\mu ^j$  space-time dependent. The combination in equation (\ref{covder}) is the well-known covariant derivative, and $B_\mu$ is the well-known gauge field. Under the gauge transformation, $B_\mu$ transforms via
\be \delta B_\mu = 2i[B_\mu , \Theta ^i T^i]-\frac{2}{g}(\partial _\mu \Theta ^i )T^i \; .\ee
With these transformations, we have 
 \be \label{inv3}\delta [(\partial _\mu -ig B_\mu)\psi] =-2i\Theta ^i T^i ((\partial _\mu -ig B_\mu)\psi) \; ,\ee
 as would be expected from equation (\ref{deltapsi}). 
 
 We  define the field strength $F_{\mu \nu}$ by 
 \bea F_{\mu \nu}&=&\partial _\nu B_\mu - \partial _\mu B_\nu + ig [B_\mu , B_\nu] \non \\
 &=&2f_{\mu \nu}^iT^i \; ,
 \eea
 where $f_{\mu \nu}^i$ are space-time dependent and  the commutator $[B_\mu , B_\nu]$ is the one defining the Lie algebra $\g$. Under the gauge transformation, $F_{\mu \nu}$ transforms by
 \be \delta F_{\mu \nu}=2i[F_{\mu \nu}, \Theta ^iT^i] \; .\ee
 Now we have all the necessary ingredients to write the Lagrangian:
 \be \label{YMlang}L=-\frac{1}{4}f_{\mu \nu}f_{\mu \nu}-\bar \psi \gamma _\mu (\partial _\mu -2ig b_\mu^iT^i)\psi -m\psi \bar \psi \; .\ee
 The Lagrangian is invariant under 
 gauge transformations. From this Lagrangian, the equations of motion of the gauge theory are derived.

\newcommand{\eab}{e_a\wedge e_b}

 \subsection{LATKe Yang-Mills, or pure Yang-Mills with matter}\label{ttYM}
 
 We can now generalize  gauge theory by replacing Lie algebras with LATKes and replacing representations of Lie algebras with representations of LATKes.

That means  we now let $\psi$ be a field in a representation $\rho$ of a LATKe as defined in equation (\ref{latkerep}), and we let $\Theta _{ab} (\eab)\in \Lambda ^2\L$ act on $\psi$ via that representation  in lieu of the action of $\Theta ^iT^i$  of the standard case of Section \ref{YMsans}. Here, $\Theta _{ab}$, $a,b=1, \ldots , \dim \L $ are space-time dependent fields antisymmetric in $a$ and $b$.  
We replace every occurrence of $\Theta ^iT^i$  in Section \ref{YMsans} by $\Theta _{ab} (\eab)$, and every index $\{ i\}$  by an antisymmetric double-index $\{ ab \}$. 

If we now inspect  the resulting equations -- which are the equations of LATKe Yang-Mills theory --  we find that the LATKe $\L$ appears  only through $\Lambda ^2 \L$. This is so because we were using representations of the LATKe in the construction, and those involve  $\Lambda ^2 \L$ rather than $\L$ (see Definition 5.1). 

We may now note that for  $\L=\L_3$, $\Lambda ^2\L$ is isomorphic to $\g _\L$ (the map $\omega$ is an isomorphism in this case) so that $\Lambda ^2\L$ is in fact the Lie algebra $\g _\L$.
 And, we also observe that  the way in which $\Lambda ^2\L$ appears in our LATKe Yang-Mills theory is  precisely the same as the way ordinary Lie algebras appear in traditional Yang-Mills theory, i.e. exactly as in Section \ref{YMsans}. 
It turns out that the definition of representations of a LATKe has conspired with the structure of the Lie algebra $\g _\L$ of the LATKe to turn LATKe Yang-Mills theory into a conventional Yang-Mills theory with Lie algebra $\g _\L$! 
And, now it is inevitable that we would think of LATKe Yang-Mills theory as a conventional Yang-Mills theory with gauge group $\exp (\g _\L)$.

Yet, there is an essential and crucial difference between conventional Yang-Mills and LATKe Yang-Mills:
in conventional Yang Mills theory, we have what is known as "pure Yang-Mills theory," where the gauge fields $B_\mu$, which live in the adjoint representation of the gauge group, are the only fields. There are no matter fields -- that is, no field $\psi$ appears -- and the Lagrangian consists only of the first term of equation (\ref{YMlang}). In general, for physical theories to include matter fields they typically have to be put in by hand. 

But in the LATKe Yang-Mills theory, this is not the case. Built into the theory is not just the adjoint representation $\Lambda ^2\L$ of $\g_\L$, but also the adjoint representation of the LATKe itself, i.e. $\L$. This representation is in fact a matter representation of $\g _\L$ and an inseparable part of {\sl pure} LATKe  Yang-Mills theory. 

Therefore, unlike pure Yang-Mills theory, pure {\sl {LATKe}} Yang-Mills theory {\sl {automatically includes matter}}, without the need to put it in by hand.

\subsection{Gauge theory for $\L_3$}
For  $\L=\L_3$, we have $\g_\L =\mathfrak{so}_4\simeq \mathfrak{sl}_2\times \mathfrak{sl}_2$  and $\L$ forms the $(2,2)$ representation (see Theorem \ref{bigtheorem}). The unitary version of  $\exp (\g _\L)$ is   $SO(4)$ or $SU(2)\times SU(2)$. The pure LATKe Yang-Mills theory for $\L _3$ is therefore an $SO(4)$ or $SU(2)\times SU(2)$ gauge theory with matter in the $(2,2)$ representation.

\subsection{LATKe Yang-Mills theory from $G_2$ manifolds} \label{YMG2}

Here we show that $\cz$, which is the singularity corresponding to the LATKe $\L _3$, indeed arises in a $G_2$ space, as we encountered in \cite{Friedmann:2002ct}.

Let $X$ be the manifold of $G_2$ holonomy which is asymptotic at infinity to a cone over $Y=\S^3\times \S ^3=SU(2)^3/SU(2)_\Delta$, where the equivalence relation $\Delta$ is $(g_1, g_2, g_3)\sim (g_1h, g_2h, g_3h)$, $g_i, h\in SU(2)$ \cite{bryant, Gibbons:1989er,Atiyah:2001qf,Friedmann:2002ct}. 

Let $\gamma \in \Z_3$ act on $Y$ as follows:
\be \gamma : (g_1, g_2, g_3) \longmapsto (\gamma g_1, \gamma g_2, g_3)\;\;  ;\; \; \;  \gamma = \left ( \begin{array} [h] {cc} e^{2\pi i/3}&0 \\ 0&e^{-2\pi i/3} \end{array}\right ) .\; 
\ee
The metric of $X$ is preserved under this action \cite{Atiyah:2001qf,Friedmann:2002ct}. Using the equivalence relation $\Delta$ to set $g_1=1$, we rewrite this action as
\be (1, g_2, g_3) \longmapsto (1, \gamma g_2\gamma ^{-1}, g_3\gamma ^{-1}) \; .
\ee
We can obtain $X$ from $Y$ by filling in one of the $SU(2)$ factors to a ball that includes the origin (recall that $SU(2)\sim \S^3$). Let us fill in the third $SU(2)$ factor, and study the singularity at the origin. 

We may write 
 \be g_2=\left ( \begin{array} [h] {cc} z_1&z_2 \\ -\bar z _2&\bar z _1 \end{array}\right )  , \; \; |z_1^2|+|z_2^2|=1\; ; \; \; \; \; \; 
         g_3=\left ( \begin{array} [h] {cc} w_1& -\bar w _2\\ w _2&\bar w _1 \end{array}\right ) \; ,
 \ee
where $z_i$ and $w_i$ are complex variables.  Then the action of $\gamma$ becomes
 \be \gamma : (z_1, z_2, w_1, w_2)\longmapsto (z_1,  e^{4\pi i/3}z_2, e^{4\pi i/3}w_1, e^{4\pi i/3}w_2)\; . \ee
 This action is singular at $z_2=w_1=w_2=0$. The locus of the singularity is the circle $z_1=e^{i\theta}$. 
 
 The singularity itself is equivalent to the one which we used in Section \ref{L3} to construct $\L_3$, with the $\C ^3$ given by coordinates $\{ z_2, w_1, w_2\}$ and $\epsilon = e^{4\pi i/3}$ (see equation (\ref{epsilonaction})).

 Therefore, we argue that the theory obtained from M-theory compactified on a $G_2$ space with a circle of $\cz$ singularities is governed by the LATKe Yang-Mills theory we constructed in Section \ref{ttYM}, i.e. it is an $SU(2)\times SU(2)$ or $SO(4)$ gauge theory on $M^4\times \S^1$ with one matter field in the $(2,2)$ representation of $SU(2)\times SU(2)$, which is the vector representation of $SO(4)$. Since this compactification is supersymmetric (it has $G_2$ holonomy), the LATKe Yang-Mills is also supersymmetric (${\cal N} = 1$). 
 
 Similarly, we argue that the same physics would result if a $\cz$ singularity appears in a Calabi-Yau space on which a string theory is compactified.
  
\subsection{Unbreakability of $\g _\L$}\label{unbroken}

We have shown that $\L_3$ leads to an $\su \times \su$ gauge theory with matter. A physicist conditioned to search for the standard model is immediately led to the following question: can we break $\su \times \su$  to $\su \times \mathfrak{u}(1)$, which is the gauge symmetry of electroweak theory and part of the standard model?

Recall that the commutators of $\L_3$ were determined from the map  $\phi : \Lambda ^2 V \rightarrow \g$, defined in equation (\ref{phi4}), which was required to satisfy the antisymmetry condition in equation (\ref{antisym123}) or (\ref{antisym}). Now we show that if $\su \times \su$ is broken to $\su \times \mathfrak{u}(1)$, {\it none} of the commutators of $\L_3$ are well-defined; in other words, they disappear. 

Let $\{ u_1, \ldots , u_6 \}$ be a basis for $\mathfrak{so}(4)$ for which the $\su \times \su$ structure is explicit:
\bea u_1=e_{12}+e_{34},&&u_2=-e_{13}+e_{24}, \hskip .8cm u_3=e_{14}+e_{23}, \\
u_4=e_{12}-e_{34},&&u_5=e_{13}+e_{24}, \hskip 1cm u_6=e_{14}-e_{23}.
\eea
Here, $\{u_1, u_2, u_3\}$ span one $\su$ factor, and $\{u_4, u_5, u_6\}$ span the other $\su$ factor. In this basis, the map $\phi$ of equation (\ref{phi4}) becomes
\bea \phi (e_1\wedge e_2)&=&(u_1-u_4)/2 \non \\ 
\phi (e_1\wedge e_3)&=&-(u_2+u_5)/2 \non \\ 
\phi (e_1\wedge e_4)&=&(u_3-u_6)/2 \non \\ 
\phi (e_2\wedge e_3)&=&(u_3+u_6)/2 \label{uphi}\\ 
\phi (e_2\wedge e_4)&=&(u_2-u_5)/2 \non \\ 
\phi (e_3\wedge e_4)&=&(u_1+u_4)/2 \non 
 \eea

Now let $u_4=u_5=0$ so that we are left with $\{u_1, u_2, u_3, u_6\}$ which forms a basis for $\su \times \mathfrak{u}(1)$. Then a fundamental requirement for the LATKe commutator is violated: 
\be  (\phi (v_i\wedge v_j))\cdot v_k\neq -(\phi (v_k\wedge v_j))\cdot v_i \hskip 1cm \forall  v_i, v_j, v_k\in V \, .\ee
For example, 
\be \phi (e_1\wedge e_2)\cdot e_3=\frac{u_1}{2} \cdot e_3=\frac{(e_{12}+e_{34})}{2} \cdot e_3=\frac{e_4}{2} \ee
while
\be \phi (e_3\wedge e_2)\cdot e_1=-\frac{(u_3+u_6)}{2} \cdot e_3=-e_{14} \cdot e_1=-e_4 \, ,\ee
so $\phi (e_1\wedge e_2)\cdot e_3\neq -\phi (e_3\wedge e_2)\cdot e_1$. The same can be checked for other combinations of $e_i$. 

So there is no longer a well-defined LATKe commutator and not even a sub-LATKe remains.

Another way of stating this result is that $\g_\L$ is unbreakable as long as $\L$ is present; or, that $\L$ protects $\g _\L$ from being broken. 
This unbreakability of $\g_{\L}$ may remind one of some global symmetries which may not be broken under certain conditions \cite{Vafa:1983tf}. 
 One may be tempted to interpret the LATKe to be a manifestation of these conditions.

 \section{Vacuum Selection Mechanism}\label{sum}
 
 \renewcommand{\baselinestretch}{1.3}
\normalsize

There was great excitement in the physics community in the 1980's when it was discovered, through a study of anomaly cancellation, that string theory came along with gauge groups -- either $E_8\times E_8$ or $SO(32)$; this discovery allowed for the hope that string theory might have some applications to phenomenology, which is governed by gauge theories \cite{Green:1984sg,Gross:1984dd}.

In the decades that followed, 
a great number of attempts at engineering a Calabi-Yau or $G_2$ space were carried out with the purpose of obtaining  a  theory in four dimensions that is as close as possible to the standard model. 
As it happened, orbifolds were employed in  Calabi-Yau compactifications of heterotic string theory for this purpose, since they induced gauge symmetry breaking by Wilson lines \cite{Hosotani:1983xw,Hosotani:1983vn,Candelas:1985en}, making the gauge group  closer to the standard model group. They also reduced the number of fermion generations that arise from the compactification, bringing that number closer to the phenomenological value of three.

Since then, it has been realized \cite{Douglas:2003um} that there is a staggering number of Calabi-Yau or $G_2$ spaces, making up what is called today the  "string landscape". Therefore, the idea of  a "vacuum selection mechanism," which is some principle  that  would  single out one vacuum or at least narrow down the choices considerably, has been sought after.

The uniqueness of the LATKe is  a vacuum selection mechanism. The selected compactification space is a Calabi-Yau or $G_2$ space with a $\cz$ singularity, and the selected vacuum theory is a supersymmetric $\su \times \su $ gauge theory with matter in the $(2,2)$ representation. 

While it has been accepted  that no vacuum selection mechanisms have as yet been proposed \cite{Douglas:2003um}, in retrospect we claim that before the present work, there did exist a  vacuum selection mechanism: anomaly cancellation. It selected a string theory with gauge group either $E_8\times E_8$ or $SO(32)$. 

While neither the uniqueness of the LATKe nor anomaly cancellation actually selects the standard model itself, our unique, simple LATKe Yang-Mills is tantalizingly close to the standard model. 

 \renewcommand{\baselinestretch}{1.4}
\normalsize

 \vskip 1cm
\no \begin{underline}{\texttt  {Acknowledgements}} \end{underline}
I am deeply indebted to I. M. Singer for inspiration and advice during the various stages of this work and for reading and commenting on an earlier draft. I am deeply indebted to F. Wilczek for discussions,  inspiration, nurturing, and mentoring. 
I wish to thank R. Jackiw, P. Deligne, V. Kac, M. Goresky, D. Kazhdan, R. MacPherson, D. Freedman,  J. Harris, M. Artin, P. Kronheimer, X. de la Ossa, J. de Jong, R. L. Jaffe, B. Zwiebach, S. Steadman, E. Farhi, W. Taylor, and C. Vafa for helpful discussions. I am very grateful to Y.P. Lee for helpful discussions and for reading and commenting on a draft. 
Last, but not least, I am grateful to A. Wiles for hospitality and encouragement. This work was supported in part by funds provided by the U.S. Department of Energy (D.O.E.) under cooperative research agreement DE-FC02-94ER40818.

\vskip 1cm
\appendix

\section{Blow-up of $\CZ$}  \label{blow}

Let $\Z_n$ be the multiplicative group generated by $\epsilon _n I_n$, where $\epsilon _n=e^{2\pi i/n}$ and $I_n$ is the $n\times n$ identity matrix. Let $\epsilon _n\in \Z_n$ act on $(z_1, \ldots ,z_n)\in \C^n$  by
\be \label{Xn} (z_1, \ldots ,z_n) \to (\epsilon _n z_1, \ldots ,\epsilon _n z_n).\ee
Denote the equivalence classes in the quotient  $X_n=\C^n/\Z_n$ by $[z_1, \ldots ,z_n]$. 

The origin of $\C^n$ is fixed under this action. The resolution of the singularity at the origin is given as follows. 

Let $Y_n=(\C^{n+1}-\{0\})/\C^*$, where $\lambda \in \C^*$ acts via
\be \label{lambda} (w_1,\ldots , w_n,w_{n+1})\to (\lambda w_1,\ldots , \lambda w_n,\lambda ^{-n}w_{n+1}) .\ee 
Denote equivalence classes in $Y_n$ by $[w_1,\ldots , w_n, w_{n+1}]$. Then $\pi:Y_n\rightarrow X_n$ is given by
\bea \label{1to1} \pi([w_1,\ldots , w_n,1])&=&[w_1, \ldots , w_n] \\
\label{exceptional}\pi([w_1,\ldots , w_n,0)&=& [0,\ldots , 0] 
\eea
Equation (\ref{1to1}) is one-to-one: the equivalence class $[w_1, \ldots , w_n,1]\in Y_n$ is determined by setting $w_{n+1}=1=\lambda ^{-n}w_{n+1}$ so now the $\lambda$ appearing in equation (\ref{lambda}) is any $n^{th}$ root of unity, leading to the same quotient action as the one defining  $X_n$ in equation (\ref{Xn}). 

Equation (\ref{exceptional}) provides us with the exceptional divisor: $\pi ^{-1}([0,\ldots ,0])$ is the set $[w_1, \ldots ,w_n,0]$, which is just $\P^{n-1}$ given the action in  equation (\ref{lambda}).

\section{The  LAnKe $\L _n$}\label{lanke}

Let $V_{n+1}$ be the standard $(n+1)$-dimensional vector representation of $\mathfrak{so}(n+1)$, and let $e_{ij}$ and $e_i$ be defined as in equations (\ref{eab}) and (\ref{ea}). Generalizing equation (\ref{phi4}), we define a map $\phi _n: \Lambda ^{n-1}V_{n+1}\rightarrow \mathfrak{so}(n+1)$ as follows:
\be \label{Ln} \phi _n (e_1\wedge \cdots \wedge \hat e_i \wedge \cdots \wedge \hat e_j \wedge \cdots \wedge e_{n+1})=(-1)^{i+j+1}e_{ij}\; , \ee
where a hat over $e_i$ means that it is omitted. 
This map yields a commutator of the $n$-th kind: 
\bea [e_1 \wedge \cdots \wedge \hat e_i \wedge \cdots \wedge e_n]&=&\phi _n (e_1\wedge \cdots \wedge \hat e_i \wedge \cdots \wedge e_{n-1})\cdot e_n\nonumber \\&=&(-1)^{i+n+1}e_{in}\cdot e_n=(-1)^{i+n+1}e_i~, \hskip .1cm i<n
\eea
\bea
[e_1 \wedge \cdots \wedge  e_{n-1}]&=&\phi _n (e_1\wedge \cdots  \wedge e_{n-2})\cdot e_{n-1} \nonumber \\ &=& (-1)^{(n-1)+n+1}e_{(n-1)n}\cdot e_n=-e_n ~.
\eea
It satisfies the requirements for a  Lie algebra of the $n$-th kind (Definition 5.2).

After a change of variables, one can show that for $\L _n$, there is a Cartan subalgebra $\mathfrak{h}_{\L_n}$ of dimension $n-1$ so $\Lambda ^{n-1}\mathfrak{h}_{\L_n}$ is one dimensional and there is a one-dimensional root space, where a root is in the dual space of $\Lambda ^{n-1}\mathfrak{h}_{\L_n}$:
\be \alpha : \Lambda ^{n-1}\mathfrak{h}_{\L_n}\longrightarrow \C \; . \ee
The  Dynkin diagram  of $\L_n$ has one node corresponding to the single cycle $\P^{n-1}$ in the exceptional divisor of the  singularity $\CZ$. 

All the mathematical definitions related to LATKes  generalize quite naturally to  LAnKes. In addition, LATKe Yang-Mills easily generalizes to LAnKe Yang-Mills, and for $\L _n$, LAnKe Yang-Mills theory is an $\mathfrak{so}(n+1)$ gauge theory with matter in the $(n+1)$-dimensional vector representation.

\section{The map $\omega$ for $B_7$} \label{omegaB3}

Let $\{E_{mn}\}_{kl}=\delta_{mk}\delta_{nl}$ be $7\times 7$ matrices. Then the following denotes a basis for $\mathfrak{so}_7$ \cite{FultonHarris,Bauerle}: 
\bea &&\epsilon_{0i}=E_{0i}-E_{i+3,0} \hskip 1.8cm (i=1,2,3) \\
&&\epsilon_{0,i+3}=E_{0,i+3}-E_{i,0} \hskip 1.3cm (i=1,2,3) \\
&&\mu_{ij}=E_{ij}-E_{j+3,i+3} \hskip 1.3cm (i,j=1,2,3)\\
&&\nu_{ij}=E_{i,j+3}-E_{j,i+3} \hskip 1.3cm (1\leq i < j\leq 3)\\
&&\rho_{ij}=E_{i+3,j}-E_{j+3,i} \hskip 1.3cm (1\leq i < j\leq 3). 
\eea
 Then $H_1=\mu_{11}$, $H_2=\mu_{22}$, $H_3=\mu_{33}$ form the Cartan subalgebra.  Recall the basis $\{ e_i, e_{j+3}, e_0 \},$ $i,j=1, 2, 3$ for $V_7$. The embedding $\psi\cdot \phi ^{-1}$ of  $\mathfrak{so}_7$ in  $C(V_7, Q)$ (Equation (\ref{embedding}))  is given by
  \bea &&\epsilon_{0i}\mapsto \frac{1}{2}e_0e_{i+3} \\
  &&\epsilon_{0,i+3}\mapsto \frac{1}{2}e_0e_i\\
  &&\mu_{ij}\mapsto \frac{1}{2}e_ie_{j+3}- \frac{1}{2}\delta_{ij}\\
  &&\nu_{ij}\mapsto \frac{1}{2}e_ie_j \\
&&\rho_{ij}\mapsto \frac{1}{2}e_{i+3}e_{j+3}~.
  \eea
 Let $f_1, \ldots, f_8$ be the basis for  $\Lambda ^\bullet W$ with $f_1=1$, $f_2=e_1\wedge e_2$, $f_3=e_1\wedge e_3$, $f_4=e_2\wedge e_3$, $f_k=e_{k-4}$ for $k=5,6,7$, and $f_8=e_1\wedge e_2\wedge e_3$. Then equation (\ref{startomega}) becomes
 \be \omega (f_1\wedge f_5)=\rho_{23}.\ee
 Using the action of $\mathfrak{so}_7$ on $\Lambda^2(\Lambda ^{\bullet}W)$ together with the intertwining condition on $\omega$ then yields
 \bea \label{omegamap} &&\omega (f_1\wedge f_5)=\rho_{23}; ~\omega (f_1\wedge f_6)=-\rho_{13}; ~\omega (f_1\wedge f_7)=\rho_{12}; \nonumber \\
 &&\omega (f_2\wedge f_5)=2\mu_{13}; ~\omega (f_2\wedge f_6)=2\mu_{23}; ~\omega (f_2\wedge f_8)=4\nu_{12}; \nonumber \\
 &&\omega (f_3\wedge f_5)=-2\mu_{12}; ~\omega (f_3\wedge f_7)=-2\mu_{32}; ~\omega (f_3\wedge f_8)=4\nu_{13}; \nonumber \\
 &&\omega (f_4\wedge f_6)=2\mu_{21}; ~\omega (f_4\wedge f_7)=2\mu_{31}; ~\omega (f_4\wedge f_8)=4\nu_{23}; \nonumber \\
  &&\omega (f_1\wedge f_2)=-\epsilon_{03}; ~\omega (f_1\wedge f_3)=\epsilon_{02}; ~\omega (f_1\wedge f_4)=-\epsilon_{01}; \nonumber \\
 &&\omega (f_2\wedge f_3)=2\epsilon_{04}; ~\omega (f_2\wedge f_4)=2\epsilon_{05}; ~\omega (f_3\wedge f_4)=2\epsilon_{06}; \nonumber \\
 &&\omega (f_5\wedge f_6)=\epsilon_{03}; ~\omega (f_5\wedge f_7)=-\epsilon_{02}; ~\omega (f_5\wedge f_8)=-2\epsilon_{04}; \nonumber \\
 &&\omega (f_6\wedge f_7)=\epsilon_{01}; ~\omega (f_6\wedge f_8)=-2\epsilon_{05}; ~\omega (f_7\wedge f_8)=-2\epsilon_{06}; \nonumber \\
 &&\omega (f_1\wedge f_8)=-\mu_{11}-\mu_{22}-\mu_{33}; ~\omega (f_2\wedge f_7)=-\mu_{11}-\mu_{22}+\mu_{33};\nonumber\\
 &&\omega (f_4\wedge f_5)=\mu_{11}-\mu_{22}-\mu_{33}; ~\omega (f_3\wedge f_6)=\mu_{11}-\mu_{22}+\mu_{33}.
 \eea
 
\vskip 1cm
\no Note added: Refs. \cite{Filipov}-\cite{takhtajan} were brought to the author’s attention after this work was completed and posted.

\end{document}